\documentclass{article}
\usepackage{authblk}
\usepackage{hyperref}
\usepackage[utf8]{inputenc}
\usepackage{fullpage}

\usepackage{caption}
\usepackage{subcaption}

\usepackage[utf8]{inputenc}
\usepackage[english]{babel}
\usepackage{csquotes}
\usepackage{hyperref}
\usepackage{url}

\usepackage{amssymb}
\usepackage{dcolumn}
\usepackage{enumerate}
\usepackage{mathdots}
\usepackage{tabularx}

\usepackage{amsmath}
\usepackage{amsthm}
\usepackage{tikz} 
\usepackage{float}
\usepackage{mathtools}
\usepackage{graphicx}

\usepackage{pgfplots}

\usepackage[ruled, commentsnumbered, longend]{algorithm2e}
\SetKwComment{Comment}{/* }{ */}

\usepackage{graphicx}
\graphicspath{ {./Images/} }

\usepackage{pgf,tikz}
\usepackage{mathrsfs}
\pgfplotsset{compat=1.18} 

\numberwithin{equation}{section}

\newtheorem{thm}{Theorem}[section]

\newtheorem{prop}[thm]{Proposition}
\newtheorem{cor}[thm]{Corollary}
\newtheorem{lem}[thm]{Lemma}

\theoremstyle{definition}
\newtheorem{defn}{Definition}[section]

\newtheorem{example}[defn]{Example}

\newcommand{\V}{\textbf{v}}
\newcommand{\one}{\textbf{1}}
\newcommand{\Grid}{\text{Grid}}

\newcommand{\len}{\text{len}}

\title{An Algorithm to Enumerate Grid Signed Permutation Classes}

\author{Sa\'ul A. Blanco\thanks{Email: \href{mailto:sblancor@indiana.edu}{sblancor@indiana.edu}}\;}
\author{Daniel E. Skora\thanks{Email: \href{mailto:danskora@iu.edu}{danskora@iu.edu}}}
\affil{Department of Computer Science\\ Indiana University\\Bloomington, IN 47408}

\date{May 31,  2023}

\begin{document}

\maketitle

\begin{abstract}
  In this paper, we present an algorithm that enumerates a certain class of signed permutations, referred to as grid signed permutation classes. In the case of permutations, the corresponding grid classes are of interest because they are equivalent to the permutation classes that can be enumerated by polynomials. Furthermore, we apply our results to genome rearrangements and establish that the number of signed permutations with fixed prefix reversal and reversal distance is given by polynomials that can be computed by our algorithm. 
\end{abstract}


\maketitle

\section{Introduction}

A permutation $\sigma$ of length $n$ \emph{contains} a permutation $\pi$ of length $k$, with $k\leq n$, if $\sigma$ has a substring $s$ where the characters of $s$ are in the same relative order as the characters of $\pi$, seen as a string. For example, $\sigma=6127354$ contains $4132$ since the characters of the substring $6254$ of $\sigma$ have the same relative order as $4132$. If $\sigma$ contains $\pi$, we write $\sigma\geq\pi$. A \emph{downset} of a permutation is the set of all permutations that are contained in $\sigma.$ A \emph{permutation class} is a set closed under containment. In other words, $\mathcal{C}$ is a permutation class if for every $\sigma\in \mathcal{C}$, $\sigma\geq\pi$ implies that $\pi\in\mathcal{C}.$ If we use $S_n$ to denote the set of all permutations of length $n$, and $\mathcal{C}$ a permutation class, studying  $|\mathcal{C}\cap S_n|$, say by providing a bound or a close expression, is a well-known combinatorial problem (see, for example~\cite{HV16,HV06,KK02, V15} and references within.) One of the most celebrated results is the proof by Markus and Tardos of the Stanley--Wilf Conjecture in~\cite{MT04} stating that unless $\mathcal{C}$ is the set of all permutations, $|\mathcal{C}\cap S_n|$ has at most exponential growth. One question that has gathered interest is which permutation classes have polynomial growth~\cite{AAB06} or those that can be enumerated by polynomials~\cite{KK02, HV06, HV16}. In answering this question, another celebrated result in the area is the so-called \emph{Fibonacci dichotomy}, first proved by Kaiser and Klazer~\cite{KK02}, stating that for every permutation class $\mathcal{C}$, then either $|\mathcal{C}\cap S_n|\geq F_n$ for all $n$, where $F_n$ denotes the $n$th Fibonacci number, or $|\mathcal{C}\cap S_n|$ is given by a polynomial for sufficiently large values of $n$. Furthermore, by combining results from~\cite{AAB06} and~\cite{HV06}, it is also true that the permutation classes enumerated by polynomials are \emph{exactly} the so-called \emph{grid permutation classes}. In fact, the following theorem is true.

\begin{thm}\cite[Theorem 1.3]{HV16}\label{thm:HVA}
Let $\mathcal{C}$ be a permutation class and $S_n$ be the symmetric group of all permutations of length $n\geq1$. Then the following statements are equivalent.
\begin{enumerate}
    \item[(a)] $|\mathcal{C}\cap S_n|$ is given by a polynomial for all sufficiently large $n$.
    \item[(b)] $|\mathcal{C}\cap S_n|<F_n$ for some $n$, where $F_n$ is the $n$th Fibonacci number ($F_1=1,F_2=2$, and $F_n=F_{n-1}+F_{n-2}$ if $n>2$.)
    \item[(c)] $\mathcal{C}$ is a \emph{grid class} of a finite set of \emph{peg permutations} (see~\cite{HV16} for these definitions).
\end{enumerate}
\end{thm}

The Homberger--Vatter algorithm~\cite{HV16}, hereinafter HVA, provides an efficient way to generate the polynomials that count the cardinality of these grid permutation classes. 

On the same token to permutation classes, one can define \emph{signed permutation classes} (see Section~\ref{sec:background} for the formal definitions) and ask which of them are enumerated by polynomials. In this case, we are interested in \emph{signed permutations}, which can be written as a string of $n$ distinct characters taken from the set $\{1,2,\ldots, n\}$ where each of the characters has a sign either positive or negative. Similarly to $S_n$, the set of all signed permutations, denoted by $B_n$, is endowed with a group structure under composition. In this context, $B_n$ is known as the \emph{hyperoctahedral group}. Indeed, $B_n$ is the group of symmetries of the hypercube, or the cross-polytope.

\subsection{Our contribution} 

In this paper, our main contributions are (i) the extension of the notion of \emph{grid class} to signed permutations. We (ii) establish that if $\mathcal{C}$ is one of these grid signed permutation classes, then $|\mathcal{C}\cap B_n|$ is given by a polynomial for all $n>1$. Moreover, we (iii) present and implement an algorithm, Algorithm~\ref{alg:one}, that computes the polynomial. An interesting application of our algorithm is the enumeration of permutation classes that are related to genome rearrangements and the burnt pancake problem. In particular, we (iv) show that the classes of signed permutations with a fixed prefix reversal and reversal distance are enumerated by polynomials and compute some of them.

The rest of the paper is organized as follows. In Section~\ref{sec:background}, we present the necessary background and definitions. In Section~\ref{sec:algorithm}, we present our algorithm to enumerate grid classes of signed permutations and the necessary auxiliary results to establish correctness and running time. Finally, in Section~\ref{sec:applications} we present applications of the algorithm to certain cases of interest that relate to the burnt pancake problem and genome rearrangement. Indeed, in bioinformatics, the prefix reversal distance and the reversal distance have biological meaning providing an idea of how ``far apart" from each other certain species are in evolutionary terms. We establish that the number of all signed permutations with prefix reversal distance and reversal distance of at most $k$ is given polynomials.

\section{Background and Notation}\label{sec:background}

We use $\mathbb{Z}^n_{\geq0}$, with $n\geq1$, to denote the set of $n$-dimensional vectors with nonnegative integer components. Similarly, we use $\mathbb{Z}^n_{>0}$ to denote the set of $n$-dimensional vectors with positive integer components. We use bold faces to denote vectors and we have reserved $\one$ to denote the vector where each component is 1 (the dimension will  be clear from context.) We denote the $i$th component of a vector $\V$ by $\V(i)$. Moreover, $S_n$ will denote the group of all permutations of $[n]:=\{1,\ldots,n\}$, the symmetric group of order $n!$.

\textbf{Signed permutations.} Given an integer $n\geq1$, let us define the set $[\pm n]:=\{\pm 1,\pm 2,\ldots,\pm n\}$. The set of all \emph{signed permutations} is the set of all bijections $\pi: [\pm n] \rightarrow [\pm n]$ satisfying the property that $\pi(-i)=-\pi(i)$ for all $1\leq i\leq n$. One usually denotes the set of all signed permutations of the set $[\pm n]$ by $B_n$ and writes $\sigma\in B_n$ as the string $\pi(1)\pi(2)\cdots \pi(n)$. We will use an underline to indicate negative entries; for example, $1 \underline{4}\hspace{0.02cm} \underline{5} 3 \underline{2}\in B_5$. We refer to the \emph{length} of a signed permutation $\pi=\pi(1)\cdots\pi(n)$ as the number of characters (positive or negative) in $\pi$ and write $\len(\pi)=n$ to denote that the length of $\pi$ is $n$, or, in other words, that $\pi\in B_n$

\textbf{Intervals, increasing and decreasing intervals, and monotone intervals.} Let $\pi=\pi(1)\cdots\pi(n)$ be a signed permutation. We say that the substring of contiguous elements $\pi(i)\pi(i+1)\cdots\pi(m)$, with $i\leq m$, is an \emph{interval} if the set $\{\pi(i),\pi(i+1),\ldots,\pi(m)\}$ equals a set of $m-i+1$ consecutive integers when we disregard their signs. If $\pi(i)\pi(i+1)\cdots\pi(m)$ is an interval, we say that $\pi(i)\pi(i+1)\cdots\pi(m)$ is \emph{increasing} if 
\[|\pi(i)|<|\pi(i+1)|<\cdots<|\pi(m)|\] and \emph{decreasing} if \[|\pi(i)|>|\pi(i+1)|>\cdots>|\pi(m)|.\]

If $\pi(i)\pi(i+1)\cdots\pi(m)$ is either increasing or decreasing, we say that $\pi(i)\pi(i+1)\cdots\pi(m)$ is a \textit{monotone interval}, or simply \emph{monotone}.

\textbf{Order isomorphic and standardization.} Given a signed permutation $\sigma\in B_n$ and $\pi\in B_k$ with $k,n\in\mathbb{Z}$ and $1\leq k\leq n$, we say that $\pi=\pi(1)\pi(2)\cdots\pi(k)$ is order isomorphic to $\sigma(i_1)\sigma(i_2)\cdots\sigma(i_k)$ if the characters of $\pi(1)\pi(2)\cdots\pi(k)$ are in the same relative order as the characters of $\sigma(i_1)\sigma(i_2)\cdots\sigma(i_k)$ and $\pi(j),\sigma(i_j)$ have the same sign for $1\leq j\leq k$. For example, $1\underline{3}2$ is order isomorphic to $2\underline{4}3$ in $12\underline{4}\hspace{0.02cm}\underline{5}3$. Furthermore, Given a string $s$ of $m$ distinct integers, each of them with a sign plus or minus, the \emph{standardization} of $s$ is the unique signed permutation $\sigma$ of length $m$ that is order isomorphic to $s$. For example, if $s=9\underline{7}43\underline{5}$, its standardization is $5\underline{4}21\underline{3}$.

We will now define a partial ordering on signed permutations, endowing $B_n$ with a poset structure.

\begin{defn}[Containment and signed permutation classes]
    We say the signed permutation $\sigma$ of length $n$ \emph{contains} the signed permutation $\pi$ of length $k\leq n$ if there exist indices $1 \leq i_1 < \cdots < i_k \leq n$ such that $\sigma(i_1)\cdots \sigma(i_k)$ is order isomorphic to $\pi=\pi(1)\cdots\pi(k)$. In particular, $\sigma(i_j)$ and $\pi(j)$ have the same sign for all $1\leq j\leq k$.
\end{defn}
If $\sigma$ contains $\pi$, then we write $\pi\leq\sigma$. Intuitively, $\pi\leq\sigma$ if we can remove entries from $\sigma$ (performing standardization as needed) to obtain $\pi$. For example, take $\sigma$ = 4\underline{1}53\underline{2} and notice that $\sigma$ contains 3\underline{1}4\underline{2}, which is obtained by removing the ``3" in $\sigma$.

 A \emph{signed permutation class} is a set of signed permutations that is closed with respect to containment. That is, if $\mathcal{C}$ is a signed permutation class with $\sigma\in \mathcal{C}$ and $\sigma\geq\pi$, then $\pi\in\mathcal{C}$.

\begin{defn}[Downsets]
Given a signed permutation $\sigma\in B_n$, we say that the \emph{downset} of $\sigma$ is the set
\[
\{\pi\in B_n\mid \pi\leq\sigma\}.
\]
\end{defn}

Another idea which is central to the behavior of our algorithm is the concept of inflation.

\begin{defn}[Inflation]
    Given a signed permutation $\pi\in B_n$ and $\V\in\mathbb{Z}_{\geq0}^n$, the process of \emph{inflating $\pi$ by $\V$} consists of replacing each character $\pi(i)$ by a monotone interval of length $\V(i)$. The interval must be increasing if $\pi(i)$ is positive and decreasing if $\pi(i)$ is negative. Moreover, each of the elements of the interval must have the same sign as $\pi(i)$. Finally, the order of these intervals must match the order of the entries of $\pi$.
\end{defn}
To inflate $\pi$ by $\V$ we write $\pi[\V]$.  For example, $\underline{1}2 [\langle 3,4\rangle] = \underline{3}\hspace{0.02cm}\underline{2}\hspace{0.02cm}\underline{1}4567$. We remark that inflating by 0 is used to remove characters; for example, notice that $21\underline{3}[\langle2,3,0\rangle]=45123$, and there is no interval corresponding to $\underline{3}$.  Further, it can be seen that $\pi \leq \sigma$ if and only if $\sigma[\V] = \pi$ where $\V$ is some vector with each component being a 0 or a 1.

\begin{defn}[Filling]
  If $\sigma=\pi[\V]$ and $\V\in\mathbb{Z}_{>0}^m$, we say that $\sigma$ \emph{fills} $\pi$. Essentially, the process of \emph{filling} a signed permutation $\pi$ is the same as that of inflating it with a vector $\V$ for which all entries are positive. In this context, we also say that $\V$ \emph{fills} $\pi$.
\end{defn} This is only a slightly stricter constraint than inflation, but it guarantees that any signed permutation obtained by filling $\pi$ will have the same general ``shape" as $\pi$.

\begin{defn}[Grid class]
  The \emph{grid class} of a signed permutation $\pi$, written $\Grid(\pi)$, is the set $\{\pi[\V] : \V\in\mathbb{Z}^n_{\geq0}\}$. Moreover, the grid class of a set of signed permutations $\Pi$, written $\Grid(\Pi)$, is the set $\{\pi[\V] : \pi\in\Pi, \V\in\mathbb{Z}^n_{\geq0}\}$.
\end{defn}

In other words, $\Grid(\Pi)$ is the set of all signed permutations that can be obtained by inflating each $\pi\in\Pi$. 
 The significance of grid classes is their relationship with signed permutation classes.  In fact, we have the following observation.
 
\begin{prop}
Let $\Pi$ be a set of signed permutations. Then, $\Grid(\Pi)$ is a signed permutation class.
\end{prop} 
\begin{proof}
    Let $\pi[\V]\in\Grid(\Pi)$. If $\tau\leq\pi[\V]$, then $\tau$ has the form $\pi[\V']$, where $\V'$ is some vector component-wise smaller than $\V$. Indeed, we can obtain $\V'$ by subtracting from each component of $\V$ the number of characters of the corresponding monotonic interval removed from $\pi[\V]$ to obtain $\tau$. For example, if $\pi=1\underline{2}3$ and $\V=\langle3,3,3\rangle$, then $\pi[\V]=123\underline{6}\hspace{0.02cm}\underline{5}\hspace{0.02cm}\underline{4}789$. If we remove the first two characters and the last four characters from $\pi[\V]$ and perform standardization, we obtain $\tau=1\underline{3}\hspace{0.02cm}\underline{2}$. Effectively, we removed from $\pi[\V]$ two elements from its first monotonic interval, one from the second one, and three from the third one. Indeed, we notice that $\tau=\pi[\langle1,2,0\rangle]$. Therefore, by definition of $\Grid(\Pi)$, $\tau=\pi[\V']$ is also in $\Grid(\Pi)$, and thus $\Grid(\Pi)$ is a signed permutation class. 
\end{proof}

\section{The Algorithm}\label{sec:algorithm}

The algorithm takes as input a finite set $\Pi$ of signed permutations. It outputs a polynomial which enumerates the elements in $\Grid(\Pi)$. In this section we describe our algorithm, explain each of its steps, discuss its running time, and illustrate the procedure with an example.

\vspace{5mm}
\noindent \textbf{\large Completion Step}

\begin{defn}[Complete sets of signed permutations]
  We say a set of signed permutations $\Pi$ is complete if every element of $\Grid(\Pi)$ can be obtained by filling some $\pi\in\Pi$. In other words, $\Pi$ is complete if for every $\sigma\in\Grid(\Pi)$ there exits a signed permutation $\pi\in \Pi$ of length $m$ and a vector $\V\in\mathbb{Z}_{>0}^m$ such that $\sigma=\pi[\V]$.
\end{defn}

We would like to insist on filling our signed permutations rather than simply inflating them because, as mentioned earlier, inflating a signed permutation by a vector with a zero entry will make it lose some of the shape.  To see this, note that any signed permutation inflated by the zero vector produces the empty permutation.  We would like to avoid double counting such permutations.

We will first attempt to replace our original set $\Pi$ with a complete set $\Pi'$ such that $\Grid(\Pi) = \Grid(\Pi')$.  Then, counting the elements in $\Grid(\Pi)$ is equivalent to counting the signed permutations that fill some element of $\Pi'$.  The following theorem tells us how to obtain this set $\Pi'$.

\begin{lem}
Any downset of signed permutations is complete.
\end{lem}
\begin{proof}
Let $\Pi$ be a downset of signed permutations and let $\sigma\in\Grid(\Pi).$ By definition of grid classes, there exits $\pi\in\Pi$ of length $m$ and a vector $\V\in\mathbb{Z}_{\geq0}^m$ such that $\sigma=\pi[\V]$. If $\pi=\pi(1)\cdots\pi(m)$, we can remove from $\pi(1)\cdots\pi(m)$ each $\pi(i)$ that corresponds to a component of $\V$ with $\V(i)=0$. After removing said entries and performing standardization, one obtains a permutation $\pi'\leq\pi$, and similarly we can obtain a vector $\V'$ by removing all the zero entries of $\V$. Since $\Pi$ is a downset, it follows that $\pi'\in\Pi$. Furthermore, $\sigma=\pi'[\V']$, so $\sigma$ fills $\pi'$.
\end{proof}

Therefore, given an input $\Pi$, the completion step of the algorithm will replace $\Pi$ by the union of the downsets of the elements of $\Pi$.

\vspace{5mm}
\noindent \textbf{\large Compacting Step}

We mentioned that the motivation for the completion step was to avoid double counting permutations obtained by inflating by a vector with one or more zero entries.  But there are some permutations we are still double counting.
Note that $12[\langle1,3\rangle]$, $12[\langle2,2\rangle]$, and $12[\langle3,1\rangle]$ all produce $1234$.

\begin{defn}[Compact signed permutations] We say that $\pi\in\Pi$ is \emph{compact} if $\pi'<\pi$ implies that $\Grid(\pi')\subset\Grid(\pi)$. We use ``$<$" and ``$\subset$" to denote strict inequality and proper subset, respectively. 
\end{defn}

For example, $12$ is not compact since $1<12$, yet $\Grid(1)=\Grid(12)$.  The definition tells us that we can safely remove such non-compact signed permutations from our set $S$ without having any effect on $\Grid(S)$.  The following proposition and lemma characterize these compact permutations and tell us what effect compactness has on double counting.

\begin{prop}\label{prop:compact}
Let $\pi$ be a signed permutation. Then the following are equivalent.
\begin{enumerate}
    \item[(a)] $\pi$ is compact,
    \item[(b)] $\pi$ does not contain an interval that is order isomorphic to $12$ or $\underline{2}\hspace{0.02cm}\underline{1}$. 
    \item[(c)] If $\pi[\V_1]=\pi[\V_2]$ where $\V_1$ has positive components and $\V_2$ has non-negative components, then $\V_1=\V_2$.
\end{enumerate}
\end{prop}

\begin{proof}
That (a) implies (b) follows from the fact that $12$ and $\underline{2}\hspace{0.02cm}\underline{1}$ are obviously not compact.  Suppose $\pi(i)\pi(i+1)$ is order isomorphic to $12$ or $\underline{2}\hspace{0.02cm}\underline{1}$.  Then any $\pi[\V]$ can be written as $\pi'[\V']$, where $\pi'$ is the result of removing $\pi(i+1)$ and $\V'$ is the result of merging $\V(i)$ and $\V(i+1)$ into one singular entry by summing them.  Since $\pi' < \pi$, $\pi$ is not compact.

Moreover, suppose (c) is true and consider the signed permutation $\sigma=\pi[\one]$, then there does not exist a signed permutation $\sigma'<\sigma$ such that $\sigma\in\Grid(\sigma')$, and therefore $\pi$ is compact. Therefore, (c) implies (a).

We are left with showing that (b) implies (c), which is proved by contraposition. Let us assume that there is a signed permutation $\pi$ such that $\pi[\V_1]=\pi[\V_2]$ where $\V_1\in\mathbb{Z}_{>0}^m$ and $\V_2\in\mathbb{Z}_{\geq0}^m$, $\V_1\neq \V_2$, and $m$ is the length of $\pi.$ Notice that $\V_1$ and $\V_2$ must differ in at least two components as otherwise the length of $\pi[\V_1]$ and $\pi[\V_2]$ would not be equal. Let $i$ be the smallest index where $\V_1(i)\neq \V_2(i)$ and $j$ be the smallest index larger than $i$ where $\V_1(j)\neq \V_2(j)$. Since $\pi[\V_1]=\pi[\V_2]$, $\V_1(i)\neq \V_2(i)$, and $\V_1(j)\neq \V_2(j)$, it follows that $\pi(i)$ and $\pi(j)$ belong to the same monotone interval in $\pi$, and since this interval is of length two at least, we know that $\pi$ contains an interval order isomorphic to either $12$ or $\underline{21}$. Therefore (b) implies (c) and the proof is now complete. 
\end{proof}

A practical implication of Proposition~\ref{prop:compact} is that condition (b) can easily be implemented to check if a permutation is compact.  All one has to do is check if $\pi(i+1)-\pi(i) = 1$ for $1\leq i\leq \len(\pi)-1$.  We make this explicit in the compacting step in Algorithm~\ref{alg:one} 

\begin{lem}\label{lem:unique}
Any signed permutation fills exactly one compact signed permutation. Namely, if $\sigma\in B_n$, there exists a unique compact $\pi\in B_m$, with $m\leq n$, and $\V\in\mathbb{Z}^m_{>0}$ satisfying $\sigma=\pi[\V]$.
\end{lem}
\begin{proof}
We will proceed by contradiction. Let us suppose that $\pi$ and $\pi'$ are two different compact permutations and $\V,\V'$ are two vectors with positive entries that satisfy $\pi[\V]=\sigma=\pi'[\V']$. Since $\pi\neq\pi'$ there exists a smallest index $i$ satisfying $\V(i)\V'(i)>1$ (which is an abbreviated way of saying that $\V(i)>1$ or $\V'(i)>1$) and $\V(i)\neq\V'(i)$. Since $\pi[\V]=\pi'[\V']$, it follows that either $\pi$ or $\pi'$ has an interval, $\pi(i)\pi(i+1)$ or $\pi'(i)\pi'(i+1)$ that is order isomorphic to $12$ or $\underline{2}\hspace{0.02cm}\underline{1}$, which contradicts the assumption that $\pi$ and $\pi'$ are compact. 
\end{proof}

Consider an arbitrary signed permutation $\pi$.  Lemma~\ref{lem:unique} tells us that $\pi$ fills exactly one compact signed permutation, and Proposition~\ref{prop:compact} (c) tells us that $\pi$ fills this compact signed permutation in a unique way.  Therefore, after the completion and compacting steps of our algorithm, we have obtained a set $S$ such that there is a bijection between members of $\Grid(S)$ and the pairs $(\pi, \V)$, where $\pi\in S$ and $\V$ fills $\pi$. In other words, we have proved the following theorem.

\begin{thm}
Let $\Pi$ be a set of signed permutations. Then
\[
\Grid(\Pi)=\bigsqcup_{\pi\in S}\Grid(\{\pi[\V]:\V\in\mathbb{Z}_{>0}^{\len(\pi)}\}),
\] where $S$ is the result of applying the completion and compacting steps of the algorithm, and ``$\sqcup$" denotes disjoint union.
\end{thm}

\vspace{5mm}
\noindent \textbf{\large Enumeration Step}

 Now, let $S$ be the set obtained after performing the completion and compacting step. To enumerate $\Grid(S)$, it suffices to enumerate the vector sets which fill each $\pi\in S$ and add them all up. For each $\pi\in S$, we wish to find the generating function for the vectors which fill $\pi$.  The generating function for vectors which fill a signed permutation of length $m$, i.e.,  all vectors  $\V\in\mathbf{Z}_{>0}^m$, is
\[\left(\frac{x}{1-x}\right)^m.\]

Therefore, we can write the generating function for $|\Grid(S)|$ as follows:
\begin{equation}\label{eq:genfunc}
    \sum_{\pi\in S} \left(\frac{x}{1-x}\right)^{\len(\pi)}.
\end{equation}

Let us recall that the generating function $\frac{x^m}{(1-x)^m}$ has a well-known explicit formula, namely:
\[
\frac{x^m}{(1-x)^m}=\sum_{i\geq0}\binom{i+m-1}{m-1}x^{i+m}.
\]

Therefore, we can easily extract the coefficient of $x^n$ in Expression~(\ref{eq:genfunc}) to obtain the enumerating polynomial for $|\Grid(\Pi)\cap B_n|$.

We include the pseudocode of our algorithm in Algorithm~\ref{alg:one}. We have also implemented our algorithm in~\cite{S22}.

\begin{algorithm}[h]
    \caption{Enumerating $\Grid(\Pi)$}    \label{alg:one}
    \KwIn{Set $\Pi$ of signed permutations}
    \KwOut{Polynomial $P(n)$ enumerating $|\Grid(\Pi)\cap B_n|$ }
   \Comment{Completion Step: Add all permutations  $\leq \pi\in S$ in the containment order.}
    Set $S =\emptyset$\;
    \For{$\pi \in \Pi$}{
        \begin{minipage}{.92\linewidth}
        Add to $S$ all permutations $\pi'\leq\pi$
        \end{minipage}
    }
    \Comment{Compacting Step: Remove all permutations $\pi\in S$ with $12,\underline{2}\hspace{0.02cm}\underline{1}\leq\pi$. This is equivalent to checking if $\pi(i+1)-\pi(i)$ equals 1 for some index $i$}
    \For{$\pi\in S$}{
        \For{$i$ from $1$ to $\len(\pi)$}{
            \uIf{$\pi(i+1) - \pi(i) = 1$}{
                \begin{minipage}{.92\linewidth}
                Remove $\pi$ from $S$
                \end{minipage}
            }
        }
    }
    \Comment{Enumeration Step: Obtain the enumerating polynomial} 
    $gen\_fcn = 0$\;
    \For{$\pi \in S$}{
        $gen\_fcn \gets gen\_fcn + \displaystyle\left(\frac{x}{1-x}\right)^{\len(\pi)}$
    }
    \Comment{poly returns the polynomial that gives the coefficient of $x^n$ in gen\_fcn} 
    $P(n)\gets poly(gen\_fcn)$\;
    \Return{$P(n)$}
\end{algorithm}

\begin{example}
Let us illustrate the running of Algorithm~\ref{alg:one} with input $\Pi=\{\underline{2}13\}$. 
\begin{description}
\item[Completion] The first step is to add to $S$ all the signed permutations $\pi\leq \underline{2}13$. At the end of the completion step,  $S=\{\varepsilon,1,\underline{1},12,\underline{1}2,\underline{2}1, \underline{2}13\}$, where we use $\varepsilon$ to denote the empty permutation, which is contained in any permutation.
\item[Compacting] We now remove the signed permutations in $S$ that are not compact. By Proposition~\ref{prop:compact}, we need to remove those permutations that contain intervals that are order isomorphic to $12$ or $\underline{2}\hspace{0.02cm}\underline{1}$. In this case, we remove 12 resulting in $S=\{\varepsilon,1,\underline{1},\underline{1}2,\underline{2}1, \underline{2}13\}$.
\item[Enumeration] The generating functions are as follows
\begin{itemize}
    \item $1$ for the empty permutation $\varepsilon$
    \item $\frac{x}{1-x}=\sum_{n\geq0}x^{n+1}$ for $1$ and $\underline{1}$,
    \item $\frac{x^2}{(1-x)^2}=\sum_{n\geq0}(n+1)x^{n+2}$ for $\underline{1}2$ and $\underline{2}1$, and
    \item $\frac{x^3}{(1-x)^3}=\sum_{n\geq0}\binom{n+2}{2}x^{n+3}$ for $\underline{2}13$.
\end{itemize} Therefore, by isolating the coefficient of $x^n$ with $n\geq1$ in the expansion of the sum 
\[
1+\frac{2x}{1-x}+\frac{2x^2}{(x-1)^2}+\frac{x^3}{(1-x)^3},
\]
we obtain the polynomial
\[
P(n)=2+2(n-1)+\binom{n-1}{2}=\frac{n^2}{2}+\frac{n}{2}+1.
\] 
\end{description}
\end{example} 

\subsection{Running time of Algorithm~\ref{alg:one}}
The running time of Algorithm~\ref{alg:one} depends on the input set $\Pi$ of signed permutations.  We condense our description of $\Pi$ to an ordered pair $(n,l)$, where $n$ is the number of permutations in $\Pi$ and $l$ is the length of the longest permutation in $\Pi$.

The completion and compacting step will add all compact permutations contained by an element of $\Pi$. Each $\pi\in\Pi$ contains at most $2^l$ sub-permutations.  Therefore, at most $n\times2^l$ operations are performed, each consisting of deleting an entry of a permutation, checking if the resulting permutation is compact, and adding it to the set if so.  Deleting an entry and checking compactness are both $O(n)$ operations, and adding to the set is constant.  So the time complexity of this step is $O(n^2 2^l)$.

The enumeration step will compute a polynomial for each $\pi\in\Pi$ and compute the sum.  At most $n\times2^l$ polynomials are computed, and each takes time $O(n)$.  So the time complexity of this step is also $O(n^2 2^l)$. Therefore, we conclude that the time complexity of Algorithm~\ref{alg:one} is $O(n^2 2^l)$.

\subsection{Differences with HVA}

In this section we compare HVA that counts elements in grid classes for permutations and Algorithm~\ref{alg:one} that counts elements in grid classes for signed permutations. It turns out that the process is streamlined in the case of signed permutations.

The first relevant difference is that in our setting, we do not need \emph{peg permutations}, which are defined in~\cite{HV16} as strings of distinct characters from 1 to $n$ in which each character is decorated with a $\bullet,+,-$. In~\cite{HV16}, they inflate and fill these peg permutations with increasing intervals if a character is decorated with a $+$, with a decreasing interval if a character is decorated with a $-$, and with an interval with 1 or 0 elements if the character is decorated with a $\bullet$. In our setting, we do not need the $\bullet$ decoration, so we use signed permutations instead of peg permutations. 

HVA can be summarized in the following five steps. Their input is a set of $P$ of peg permutations.

\begin{description} 
\item[Completion] In this step, the authors add to $S$ all the peg permutations that are contained in peg permutations from $S$. Our completion step is similar but we use signed permutations instead of peg permutations.
\item[Compacting] In the compacting step, HVA removes the non-compact peg permutations from $S$. Our compacting step is similar, but the condition to check which permutations to remove is easier in Algorithm~\ref{alg:one}.
\item[Cleaning] There are some redundancies in the enumeration process in HVA that are not removed by the compacting step. The authors need to remove the peg permutations that are not \emph{clean} (see~\cite[Section 2.3]{HV16}) from $S$. This step is not needed in Algorithm~\ref{alg:one} due to Lemma~\ref{lem:unique}.
\item[Combining] In the combining step, HVA  takes the union of certain convex sets produced during cleaning. Algorithm~\ref{alg:one} does not use the convex sets used in HVA, so this step can safely be skipped for signed permutations. 
\item[Enumeration] This step is the same for both HVA and Algorithm~\ref{alg:one}.
\end{description}

Moreover, it is worth pointing out that HVA works for ``sufficiently large values of $n$." The reason why HVA requires this condition on $n$ is because the \emph{minimal filling vector} $\V$ of a peg permutation $\tilde{\rho}$ requires $\V(i)=2$ if $\tilde{\rho}(i)$ is decorated with a $\bullet$. By contrast, our minimal filling vector is always $\mathbf{1}$ and therefore our polynomials work for all values of $n$.

\section{Applications}\label{sec:applications}

\subsection{Pancake Problem}

 Given a permutation $\pi\in S_n$ written in one-line notation, the \emph{prefix reversal} $f_i$ with $2\leq i\leq n$ is a a function on permutations that ``reverses" the first $i$ characters of $\pi$. More precisely, if $\pi=\pi(1)\pi(2)\cdots\pi(n)$ then $f_i(\pi)=\pi(i)\pi(i-1)\cdots \pi(1)\pi(i+1)\cdots\pi(n)$, with $2\leq i\leq n$. In its original formulation in~\cite{Dweighter75}, the \emph{pancake problem}, asks for the minimum number of applications of prefix reversals needed to sort a given permutation (i.e., the \emph{prefix reversal distance}). In the only academic paper that Microsoft's Bill Gates ever wrote, Gates and Papadimitriou~\cite{GatesPapa} provided the first non-trivial bound for the prefix reversal distance. This bound was later improved in~\cite{Chitt}. Furthermore, the pancake problem was later extended in~\cite{GatesPapa} to the setting of signed permutations, where the prefix reversals also swap the signs of the characters from ``$+$" to ``$-$" and vice versa. In the case of $B_n$, there are $n$ prefix reversals $f^B_i$, $1\leq i\leq n$, that will reverse the first $i$ characters of a permutation written in one-line notation and switch their sign. For example, $f^B_1(13\underline{2})=\underline{1}3\underline{2}$ and $f^B_2(13\underline{2})=\underline{3}\hspace{0.02cm}\underline{1}\hspace{0.02cm}\underline{2}$. Finding the minimum number of signed prefix reversals needed to sort a signed permutation is known as the \emph{burnt pancake problem}.

It is worth mentioning what brought our attention to this project. In~\cite{BBP19Perm}, utilizing entirely enumerative methods like the principle of inclusion-exclusion and the classification of small cycles in the pancake and burnt pancake graph from~\cite{BBP19,BBP19Perm}, the authors provided a polynomial that counts the number of stack of $n$ pancakes that require exactly four flips to be sorted, as well as the number of stacks of burnt pancakes that require four flips to be sorted. Furthermore, the authors also conjectured other polynomials that would yield the number of stacks of pancakes or burnt pancakes that require exactly $k$ flips to be sorted for small values of $k$ (see~\cite[Section 6]{BBP19Perm}.) Eventually, V.~Vatter pointed out their HVA~\cite{HV16} to the authors in~\cite{BBP19Perm}. Immediately, two questions became evident.
\begin{enumerate}
    \item[Q1:] Is it possible to modify HVA to work for signed permutations?, and
    \item[Q2:] HVA missed some of the numbers obtained in~\cite{BBP19Perm} due to HVA's nature of working for sufficiently large values of $n$. Why does HVA work for sufficiently large values of $n$ and not all values of $n$?
\end{enumerate}

The answer to Q1 is ``yes" and Algorithm~\ref{alg:one} is a modified version of HVA tailored to enumerate signed permutation classes. Moreover, relating to Q2 and unlike HVA, Algorithm~\ref{alg:one} works for all values of $n$ and not just for ``sufficiently large values of $n$" due to the minimum filling vector in the case of signed permutations being always $\mathbf{1}$. In what follows of this subsection, we apply Algorithm~\ref{alg:one} to obtain the number of stacks of burnt pancakes that take at most $k$ flips to be sorted. In particular, we establish ~\cite[Conjecture 6.2, 6.3]{BBP19Perm}.

\subsubsection{From Burnt Pancake Stacks to Grid Classes}\label{sec:burnt}

We know that we can represent a stack of burnt pancakes as a signed permutation, so our next objective is to represent the stacks of burnt pancakes which can be sorted in $\leq k$ flips as a grid class.  It is easy to see that $\Grid(1) = \{\varepsilon, 1, 12, 123, ...\}$ represents the stacks of pancakes which can be sorted in $\leq 0$ flips, but what about for other values of $k$?  In this section, we will show how to choose $\Pi_k$ such that $\Grid(\Pi_k)$ is the class of signed permutations that require $\leq k$ prefix reversals to be sorted.

First, let's take a moment to consider what a grid class intuitively represents.  In this algorithm, and specifically in regard to the pancake problem, we can visualize signed permutations in two different ways, depending on context.  Suppose we are presented with a grid class $\Grid(\Pi)$. Then $\Pi$ and $\Grid(\Pi)$ are both sets of signed permutations, but we can interpret their elements differently.  If $\pi\in\Grid(\Pi)$, we think of each entry of $\pi$ as representing a burnt pancake.  Moreover, if $\pi\in\Pi$, we can also think of each entry of $\pi$ as a block of indefinitely many consecutive pancakes which is in sorted order internally.  Similar to a single burnt pancake, each block has orientation and relative order and can be treated as a singular entity.  Conceptually, when we inflate an element of $\Pi$, we are simply assigning a nonnegative integer size to each block.

Consider $\Grid(\underline{1}\hspace{0.02cm}2)$ as an example.  We interpret $\underline{1}\hspace{0.02cm}2$ as an upside down block of pancakes on top of a larger, right side up block of pancakes, and every signed permutation in $\Grid(\underline{1}2)$ fits this characterization.  For example, $\underline{1}\hspace{0.02cm}2[\langle 3,3\rangle] = \underline{3}\hspace{0.02cm}\underline{2}\hspace{0.02cm}\underline{1}456$, and we can sort $\underline{3}\hspace{0.02cm}\underline{2}\hspace{0.02cm}\underline{1}456$ similarly to how we would sort $\underline{1}\hspace{0.02cm}2$.
\begin{figure}[t]
\centering
\includegraphics{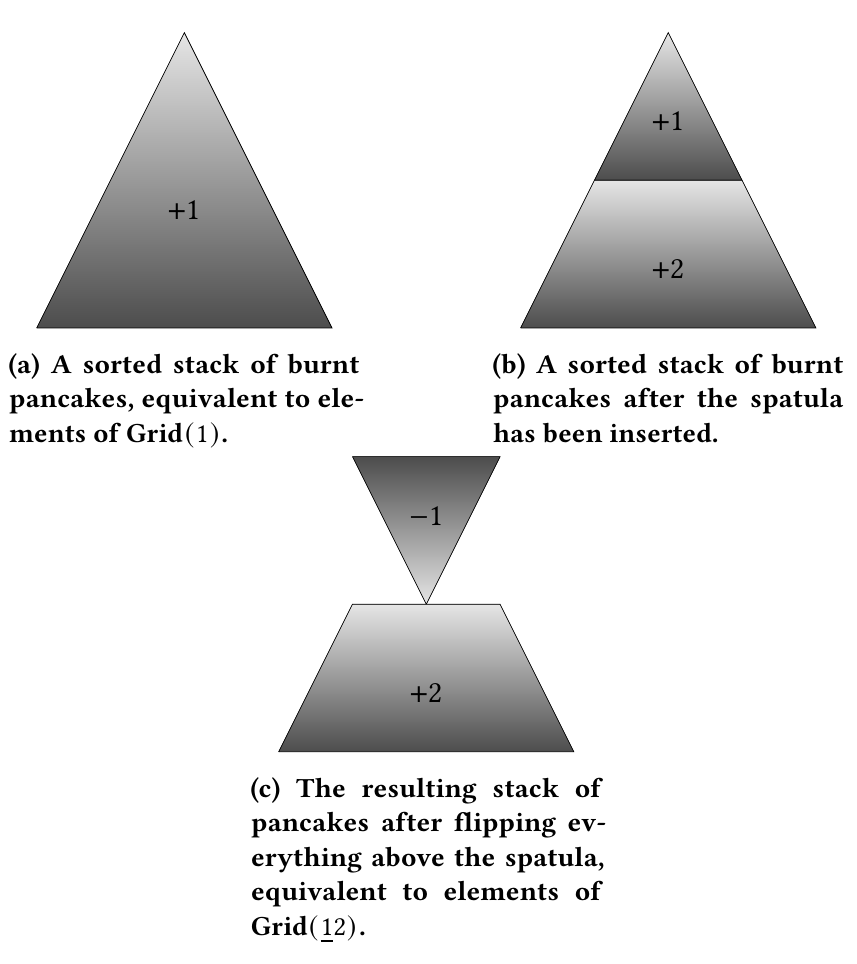}
\caption{}
\label{fig1}
\end{figure}

\begin{figure}[b]
\centering
\includegraphics{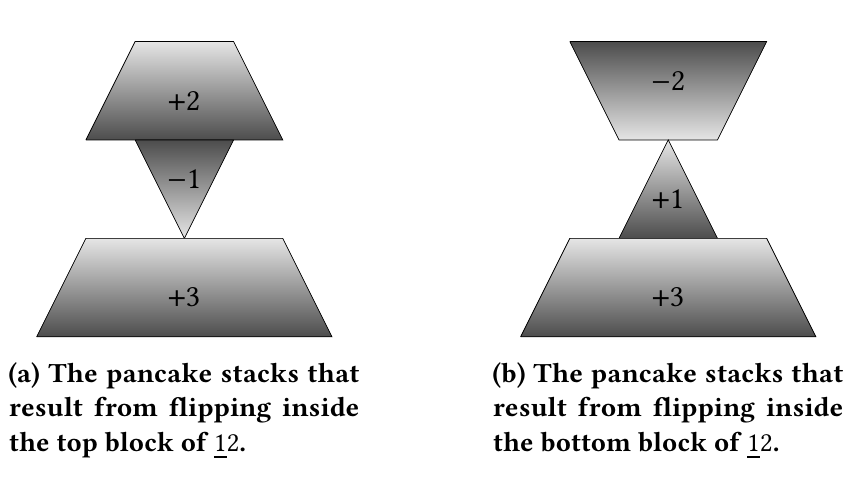}
\caption{}
\label{fig2}
\end{figure}

For the remainder of this discussion, we visualize our signed permutations as stacks of blocks rather than stacks of pancakes.  This intuition will help us choose $\Pi_k$.

To get any stack of pancakes which is one flip from being sorted, we start with a sorted stack, insert our spatula somewhere in the middle of the stack, and flip.  In other words, we start with the signed permutation $1$ from Figure~\ref{fig1}a.  Then we insert our spatula, splitting the block into two, shown in Figure~\ref{fig1}b.  Finally, we flip, which produces $\underline{1}2$ as seen in Figure~\ref{fig1}c.  So intuitively we expect $\Grid(\underline{1}2)$ to represent the stacks of pancakes which can be sorted in $\leq 1$ flips.

We can continue this process recursively.  To get any pancake stack which is two flips from the identity, we start with a stack which is one flip from the identity, insert our spatula anywhere, and flip.  But the stacks which are one flip from the identity are represented by $\underline{1}2$, which consists of two blocks.  We can flip in the middle of either block, which means we are going to obtain two different signed permutations.  If we flip inside the $\underline{1}$ block, we get $2\underline{1}3$.  If we flip inside the $2$ block, we get $\underline{2}13$.  These stacks are shown in Figure~\ref{fig2}a and Figure~\ref{fig2}b, respectively.  So we expect $\Grid(\{2\underline{1}3, \underline{2}13\})$ to represent the stacks of pancakes which can be sorted in $\leq 2$ flips.

This process is straightforward, and we can continue this recursion to generate $\Pi_k$ for any $k$.  A simple induction argument shows that $\Pi$ will contain $k!$ elements, each of length $k+1$.

The recursive procedure we have described is meant to mirror all possible ways we could rearrange a sorted stack of pancakes in $k$ flips.  However, our procedure requires us to perform exactly $k$ flips, and each flip must be in the middle of a block, splitting it into two.  Consequently, there are certain stacks of pancakes which are not obviously contained in $\Grid(\Pi_k)$: any stack generated by strictly $<k$ flips or any stack generated by at least one flip that didn't ``separate" pancakes.

Although our process does not allow us to flip at the edge of a block, it allows us to flip infinitesimally close to the edge of a block.  In other words, if we flip in the middle of a block and then inflate one of the resulting two blocks by zero, we have simulated flipping at the edge of the original block.  We can use this trick to flip in between blocks, which is equivalent to a flip that doesn't separate pancakes.  We can also use this trick to flip at the very top or very bottom of the stack, and the former is equivalent to not flipping at all.  This analysis shows us that $\Grid(\Pi_k)$ must contain all stacks of pancakes with prefix reversal distance $\leq k$, and we formalize this in a proposition.

\begin{prop}\label{prop:gridpancake}
For each $k\geq0$, the set of all burnt pancake stacks with  prefix reversal distance $\leq k$ is equivalent to a signed grid class.
\end{prop}

\begin{proof}
Specifically, this grid class is $\Grid(\Pi_k)$, which we generate recursively starting with $\Pi_0 = \{1\}$.  For each $\pi\in\Pi_k$, and for each $1\leq i \leq k+1$, we generate a new permutation $\pi_i$ by first inflating $\pi$ by the vector $e_i + \one$, and then applying the prefix reversal $f_i$ to the result.  Here $e_i$ denotes the standard basis vector where the $i$th component is $1$ and all others are $0$. Define $\Pi_{k+1}$ as the set of all these $\pi_i$.

First, we show that $\Grid(\Pi_k)$ contains every signed permutation whose prefix reversal distance is $\leq k$.  If $\pi$ has prefix reversal distance $\leq k$, then it is an inflation of a (not necessarily unique) compact permutation with prefix reversal distance $k$.  According to the discussion above, this compact permutation is contained in $\Pi_k$, so $\pi$ is in $\Grid(\Pi_k)$.

Before proving the converse, we extend the notion of \emph{prefix reversals} to apply to vectors. Applying $f_i$ to $\V$ reverses the first $i$ components of $\V$.

Now, we show that every signed permutation in $\Grid(\Pi_k)$ has a prefix reversal distance $\leq k$.  Any permutation $\sigma\in\Grid(\Pi_k)$ can be written as $\pi[\V]$ for some $\pi\in\Pi_k$ and $\V\in\mathbb{Z}^{n}_{\geq0}$.  Due to how we constructed $\Pi_k$, it is clear that $\pi$ is the result of performing a sequence $s$ of $\leq k$ prefix reversals on the identity permutation of length $k+1$.  Therefore, we can generate $\sigma$ by an analogous sequence $s'$ of prefix reversals on the identity permutation of length $\len(\sigma)$, and consequently $\sigma$ has prefix reversal distance $\leq k$.  This sequence $s'$ can be explicitly described as follows.  We traverse $s$, and for every prefix reversal $f_i$ , the corresponding prefix reversal in $s'$ is $f_x$, with $x$ being the sum of the first $i$ components in $\V$.  Then, after applying $f_x$, we apply $f_i$ to $\V$ for the sake of organizing the sizes of each inflated interval, and repeat the process with the next prefix reversal in $s$. For example, consider $\pi=\underline{2}13$ and $\V=\langle 1,2,3\rangle$. Notice that $f_1f_2(\pi)=123$ and the corresponding sequence of flips to sort $\pi[\V]=\underline{3}12456$ starts with $f_{1+2}=f_3$. Then we apply $f_2$ to $\V$ obtaining $\langle 2,1,3\rangle$, and then the corresponding prefix reversal to $f_1$ is $f_2$ since $2$ is the first component of $\langle 2,1,3\rangle$. Thus applying $f_2$ and then $f_1$ to $\pi$ corresponds to applying $f_3$ and then $f_2$ to $\pi[\V]$. Indeed, $f_2f_3(\underline{3}12456)=123456$.
\end{proof}

Now that we have established that the set of permutations with fixed prefix reversal distance is a grid class, we can apply Algorithm~\ref{alg:one} to compute some of the polynomials that enumerate them for certain values of $n$ in the section that follows. We remark that we recover~\cite[Theorem 4.3]{BBP19Perm} previously established utilizing basic counting techniques and not symbolic computation. 

\subsubsection{Results}

We used the algorithm to acquire the enumerating polynomial $R^{B}_{\leq k}(n)$, which outputs the number of distinct stacks of $n$ pancakes which take $k$ or less flips to be sorted.  We obtained polynomials for values of $k$ from 1 to 10.  Our results are below. To make things easier to display, we denote the polynomials by an array of their coefficients where the entry at index $i$, with $1\leq i\leq k+1$, represents the coefficient of $n^{i-1}$ in $R^B_{\leq k}$. For example, $R^{B}_{\leq 3}(n) =  1+n-n^2+n^3$ and we represent it by $[1,1,-1,1]$.

\begin{align*}
R^{B}_{\leq 1}(n) &= \left[1,1\right],\\
R^{B}_{\leq 2}(n) &= \left[1,0,1\right],\\
R^{B}_{\leq 3}(n) &= \left[1,1,-1,1\right],\\
R^{B}_{\leq 4}(n) &= \left[ 1,-\frac{1}{2},3,-\frac{5}{2},1\right],\\
R^{B}_{\leq  5}(n) &= \left[1,\frac{1}{2},-\frac{25}{6},\frac{17}{2},-\frac{29}{6},1\right],\\
R^{B}_{\leq 6}(n) &= \left[1,\frac{299}{30},-5,-\frac{73}{4},21,-\frac{463}{60},1\right],\\
R^{B}_{\leq 7}(n) &= \left[1,-\frac{3529}{30},\frac{24697}{120},-\frac{3167}{48},-\frac{889}{16},\frac{3569}{80},-\frac{2699}{240},1\right],\\
R^{B}_{\leq 8}(n) &= \left[1,\frac{92843}{84},-\frac{48217}{20},\frac{1230329}{720},- \frac{7787}{24},- \frac{2659}{18},\frac{10117}{120},- \frac{77323}{5040},1\right],\\
R^{B}_{\leq 9}(n) &=
\left[1, -\frac{1713461}{168}, \frac{28102741}{1120}, -\frac{3620111}{160}, \frac{52327853}{5760}, -\frac{13571}{12}, -\frac{997679}{2880}, \frac{163277}{1120}, -\frac{806941}{40320}, 1\right],\text{ and}\\
R^{B}_{\leq 10}(n) &=
\left[1, \frac{29555642}{315}, -\frac{1264975307}{5040}, \frac{11803588051}{45360}, -\frac{77767535}{576}, \frac{307180691}{8640}, -\frac{4420823}{1440}, -\frac{22399579}{30240}, \right.\\&\left.\frac{948575}{4032}, -\frac{4576633}{181440}, 1\right].
\end{align*}

If one wishes to consider the stacks of burnt pancakes with prefix reversal distance of \emph{exactly} $k$, the enumerating polynomial is $R^{B}_{\leq k}(n)-R^{B}_{\leq k-1}(n)$. For example,

\[
R^{B}_{\leq 4}(n)-R^{B}_{\leq 3}(n)=\frac{1}{2}(n-1)^2(2n-3),
\]
which is~\cite[Theorem 4.3]{BBP19Perm}. In fact, we are able to confirm that the polynomials given in~\cite[Conjecture 6.1]{BBP19Perm} are all correct. For the sake of completeness, we include the conjecture here, now a proposition. We use $R^B_{k}$ to denote the polynomial that enumerates how many signed permutations have prefix reversal distance exactly $k$.  

\begin{prop}\label{con:bn} If $n\geq1$, then
\begin{enumerate}
    \item[(i)] $R^B_5(n)=\frac{1}{6}n(n-1)(n-2)(6n^2-17n+3)$,
    \item[(ii)] $R^B_6(n)=\frac{1}{60}n(n-1)(n-2)(60n^3-343n^2+401n+284)$,
    \item[(iii)] $R^B_7(n)=\frac{1}{240}n(n-1)(n-2)(n-3)(240n^3-1499n^2+925n+5104)$,
    \item[(iv)] $R^B_8(n)=\frac{1}{5040}n(n-1)(n-2)(n-3)(5040n^4-52123n^3+113415n^2+314716n-1027242)$, and
    \item[(v)] $R^B_9(n)=\frac{1}{40320}(n-1) (n-2) (n-3) (n-4) (40320 n^5-444061 n^4+644746 n^3+6638777 n^2-18991470 n).$
\end{enumerate}
\end{prop}

Moreover,~\cite[Conjecture 6.3]{BBP19Perm} is now a corollary that follows by using the \emph{Gregory-Newton interpolation formula} for integer-valued polynomials.

\begin{cor}
 If $k,n\geq1$, then
 \begin{displaymath}R_k^B(n) =\sum_{j=1}^k \left(\sum_{i=0}^{k-j} (-1)^i \binom{i+j-1}{i}\binom{n}{i+j}\right)R_k^B(j), \text{ for } k \geq 1.\end{displaymath}
\end{cor}

\subsection{Genome Rearrangement}

A genome is the entire set of DNA instructions found in a cell. A genome contains the genetic information of an organism; for example, in humans, the genome consists of 23 pairs of chromosomes. 
The DNA is partitioned into sections called ``genes" where each gene is a genetic recipe for some very specific feature of the organism.  Assuming each gene is an asymmetrical sequence of nucleotides, we can think of the genome as a sequence of genes where each gene has orientation.

Evolution occurs when a mutation alters an organism's genome, acting on a sequence of genes.  By far the most common of these is a \emph{block reversal}, or simply \emph{reversal}, in which case an entire section of genes is reversed. More formally, if $\pi=\pi(1)\pi(2)\cdots\pi(n)\in B_n$, a block reversal $b_{i,j}$ with $1\leq i\leq j\leq n$ transforms $\pi$ into the permutation 
\[
\pi(1)\cdots\pi(i-1)\pi'(j)\pi'(j-1)\cdots\pi'(i+1)\pi'(i)\pi(j+1)\cdots\pi(n),
\] where $\pi'_k$ denotes $-\pi_k$ for $i\leq k\leq j$.

If we want to compare an organism to its ancestor, we can represent their genomes by signed permutations.  The descendant's DNA can be written as the identity permutation, and based on this labelling of genes we can write the ancestor's DNA similarly.  It is of interest to biologists to determine the minimum number of reversals which would transform a given signed permutation into the identity, since this probably represents the number of mutations which occurred during evolution. In this context, sorting permutations by reversals has a biological interpretation, which has been explored in~\cite{CohenBlum, HannenPev, H08}.

We are interested in the class of signed permutations that take $\leq k$ block reversals to be sorted (i.e., \emph{reversal distance} of $k$).  These classes are also equivalent to grid classes, which we establish by utilizing a method similar to what we did with the pancake problem in the previous section. Indeed, in the Section~\ref{sec:burnt}, we obtained sets $\Pi_k$ recursively by simulating a prefix reversal at every entry of every signed permutation in $\Pi_{k-1}$.  For this problem, we simulate block reversals by choosing \emph{two} blocks (possibly not distinct) at which to split and invert the middle portion.

\begin{prop}
For each $k\geq0$, the set of all signed permutations with reversal distance $\leq k$ is equivalent to a grid class.
\end{prop}

\begin{proof}
The grid class is $\Grid(\Pi_k)$, which we generate recursively starting with $\Pi_0 = \{1\}$.  For each $\pi\in\Pi_k$, and for each $1\leq i \leq j \leq k$, we generate a new permutation $\pi_{ij}$ by first inflating $\pi$ by the vector $e_i + e_j + \one$, and then applying the block reversal $b_{i+1,j+1}$ to the result.  Define $\Pi_{k+1}$ as the set of all these $\pi_{ij}$.

We first establish that $\Grid(\Pi_k)$ contains every signed permutation whose reversal distance is $\leq k$.  If $\pi$ has reversal distance $\leq k$, then it is an inflation of some compact permutation with reversal distance $k$. This compact permutation is contained in $\Pi_k$, so $\pi$ is in $\Grid(\Pi_k)$.

We extend the notion of \emph{block reversals} to apply to vectors.  So applying $b_{i,j}$  to $\V$ reverses the $i$th through $j$th components of $\V$, inclusively.

We now show that every signed permutation in $\Grid(\Pi_k)$ has a reversal distance $\leq k$.  Any permutation $\sigma\in\Grid(\Pi_k)$ can be written as $\pi[\V]$ for some $\pi\in\Pi_k$ and $\V\in\mathbb{Z}^{n}_{\geq0}$.  We generated $\pi$ by a sequence $s$ of $\leq k$ block reversals on the identity permutation.  Therefore, we can generate $\sigma$ by an analogous sequence $s'$ of block reversals on the identity permutation of length $\len(\sigma)$, and consequently $\sigma$ has reversal distance $\leq k$.  This sequence $s'$ can be explicitly described as follows: We traverse $s$, and for every block reversal $b_{i,j}$, the corresponding reversal in $s'$ is $b_{x,y}$, with $x$ being the sum of the first $i-1$ components in $\V$ plus $1$, and with y being the sum of the first $j$ components in $\V$.  We then apply $b_{i,j}$ to $\V$ and move on to the next block reversal in $s$.  After completing the process, we are left with the uniquely determined sequence $s'$ that  sorts $\sigma$.
%
%
\end{proof}

Now that we have established that the set of permutations with fixed reversal distance is a grid class, we can apply Algorithm~\ref{alg:one} to compute some of the polynomials that enumerate them for certain values of $n$ in the section that follows.  

\subsubsection{Results}
We used the algorithm to acquire the enumerating polynomial $P^{B}_{\leq k}(n)$, which outputs the number of elements of $B_n$ which are $\leq k$ block reversals from the identity.  We obtained polynomials for values of $k$ from 1 to 5.  Our results are below. Once again, we denote the polynomials by an array of their coefficients where the entry at index $i$, with $1\leq i\leq 2k+1$, represents the coefficient of $n^{i-1}$ in $P^B_{\leq k}$.

\begin{align*}
P^{B}_{\leq 1}(n) &= \left[1,\frac{1}{2},\frac{1}{2}\right],\\
 P^{B}_{\leq 2}(n) &= \left[1,\frac{1}{3}, \frac{1}{3}, \frac{1}{6}, \frac{1}{6}\right],\\
 P^{B}_{\leq 3}(n) &= \left[1,\frac{1}{3},\frac{35}{72},\frac{7}{48},-\frac{5}{144},\frac{1}{48},\frac{7}{144}\right],\\
 P^{B}_{\leq 4}(n) &= \left[1,\frac{131}{420},\frac{617}{1260}, -\frac{1}{120},\frac{67}{1440}, \frac{53}{240}, -\frac{17}{360}, -\frac{41}{1680}, \frac{37}{3360}\right],\text{ and}\\
      P^{B}_{\leq 5}(n) &= \left[1,\frac{331}{2520}, \frac{24727}{50400}, \frac{4703}{22680}, \frac{16945}{72576}, \frac{931}{17280}, -\frac{20059}{86400}, \frac{7267}{60480},\frac{145}{24192},-\frac{925}{72576},\frac{3767}{1814400} \right].
\end{align*}

\section{Concluding Remarks}

Our main contribution is Algorithm~\ref{alg:one} that enumerates grid signed permutation classes, which we implement in~\cite{S22}. Studying the size of permutation classes has received much attention including profound results such as the proof of the Stanley--Wilf Conjecture and the Fibonacci dichotomy. Moreover, Homberger and Vatter provided an algorithm that outputs a polynomial of degree $n$ for all the permutation classes $\mathcal{C}$ satisfying $|\mathcal{C}\cap S_n|$ being a polynomial for sufficiently large values of $n$. We extend the notion of grid permutation class to signed permutations and prove that if $\mathcal{C}$ is a grid signed permutation, then $|\mathcal{C}\cap B_n|$ is a polynomial computed by Algorithm~\ref{alg:one}. In this case, one does not need $n$ to be sufficiently large. Our grid signed permutation classes include classes that are of interest in bioinformatics, in particular to genome rearrangements.  Indeed, in nature, one way in which species evolve is by applying a block-reversal to their genome sequences. In this bioinformatics context, we obtain bounds for how many potential species are within ``distance" $k$ of one another. As future work, we expect to obtain a result similar to Theorem~\ref{thm:HVA} that applies to signed permutations.

\section*{Acknoledgments}
We thank the anonymous reviewers for their valuable comments that improved the quality of the paper. We also thank V.~Vatter who pointed out~\cite{HV16} to the first author.


\bibliographystyle{plain}

\end{document}